\documentclass{amsart}
\usepackage{amsfonts,amssymb,amsmath,amsthm,mathrsfs}
\usepackage{url}
\usepackage[dvips]{epsfig}
\newtheorem{thrm}{Theorem}[section]
\newtheorem{lem}[thrm]{Lemma}
\newtheorem{prop}[thrm]{Proposition}
\newtheorem{cor}[thrm]{Corollary}
\newtheorem{rem}{Remark}
\theoremstyle{definition}
\newtheorem{definition}[thrm]{Definition}

\begin{document}
\author{Carlo Alberto Mantica and Luca Guido Molinari} 
\address{Physics Department, Universit\`a degli Studi di Milano and I.N.F.N.,
Via Celoria 16, 20133, Milano, Italy.\\ 
Present address of C.~A.~Mantica: I.I.S. Lagrange, 
Via L. Modignani 65, 20161, Milano, Italy}
\email{luca.molinari@mi.infn.it, carloalberto.mantica@libero.it }
\subjclass[2010]{53B20, 53B21}.
\keywords{Weakly-Ricci Symmetric manifolds, Pseudo- Projective Ricci symmetric,
conformal curvature tensor, quasi conformal curvature tensor, 
conformally symmetric, conformally recurrent, Riemannian manifolds. 
Weakly Z symmetric Manifolds.}
%%%%%%%%%%%%%%%%%%%%%%%%

%
\title{Weakly Z symmetric manifolds}
\begin{abstract}
We introduce a new kind of Riemannian manifold 
that includes weakly-, pseudo- and pseudo projective- Ricci symmetric 
manifolds. The manifold is defined through a generalization 
of the so called $Z$ tensor; it is named 
{\em weakly $Z$ symmetric} and denoted by $(WZS)_n$. 
If the $Z$ tensor is singular we
give conditions for the existence of a proper concircular vector. 
For non singular $Z$ tensor, we study 
the closedness property of the associated 
covectors and give sufficient conditions for the existence of a 
proper concircular vector in the conformally harmonic case, and
the general form of the Ricci tensor.
For conformally flat $(WZS)_n$ manifolds, we derive the
local form of the metric tensor.\\
Date: 20 february 2011.
\end{abstract}
\maketitle

\section{Introduction}
In 1993 Tamassy and Binh \cite{[25]} introduced and studied a Riemannian 
manifold whose Ricci tensor\footnote{Here we define 
the Ricci tensor as $R_{kl}= -R_{mkl}{}^m $ %\cite{[29]} 
and the scalar curvature as $R = g^{ij} R_{ij} $. $\nabla_k$ is the covariant 
derivative with reference to the metric $g_{kl}$. We also put $\|\eta\|=
\sqrt{\eta^k\eta_k}.$} 
satisfies the equation:
\begin{equation}	
\nabla_k R_{jl}  = A_k R_{jl}  + B_j R_{kl}  + D_l R_{kj} .\label{eq1.1}
\end{equation}
The manifold is called {\em weakly Ricci symmetric} and denoted by
$(WRS)_n$. The covectors $A_k$, $B_k$ and $D_k$ are the 
{\em associated 1-forms}. The same manifold with the 1-form $A_k$ 
replaced by $2A_k$ was studied by Chaki and Koley \cite{[6]}, and called 
{\em generalized pseudo Ricci symmetric}. The two structures extend 
{\em pseudo Ricci symmetric} manifolds, $(PRS)_n$, introduced by Chaki 
\cite{[4]}, where $\nabla_k R_{jl}  = 2A_kR_{jl}  + A_j R_{kl}  + A_l R_{kj}$
(this definition differs from that of R.~Deszcz \cite{[14]}).\\
Later on, other authors studied the  
manifolds \cite{[9],[18],[11]}; in \cite{[11]} some global properties 
of $(WRS)_n$ were obtained, and the form of the Ricci tensor was found. 
In \cite{[9]} %conformally flat 
generalized pseudo Ricci symmetric manifolds were considered, where the 
conformal curvature tensor 
\begin{eqnarray}
C_{jkl}{}^m  = R_{jkl}{}^m  + \frac{1}{n - 2}(\delta_j{}^m R_{kl}  
- \delta_k{}^m R_{jl}  + R_j{}^m g_{kl}  - R_k{}^m g_{jl} )\nonumber\\ 
- \frac{R}{(n - 1)(n - 2)} (\delta_j{}^m g_{kl}  - \delta_k{}^m g_{jl} )
\label{eq1.2}
\end{eqnarray}
vanishes (for $n=3$: $C_{jkl}{}^m  = 0$ holds identically, 
\cite{[22]}) and the existence of a proper concircular vector was proven.
In \cite{[18]} a quasi conformally flat $(WRS)_n$ was studied, and again the 
existence of a proper concircular vector was proven.\\ 
In \cite{[2]} $(PRS)_n$ with harmonic curvature 
tensor (i.e. $\nabla _m R_{jkl}{}^m  = 0$) or with harmonic conformal curvature 
tensor (i.e. $\nabla _m C_{jkl}{}^m  = 0$) were considered.\\
Chaki and Saha considered the projective Ricci tensor $P_{kl}$, obtained by a 
contraction of the projective curvature tensor $P_{jkl}{}^m$ \cite{[15]}:
\begin{equation}
P_{kl}  = \frac{n}{n - 1}\left ( R_{kl}  - \frac{R}{n} g_{kl} \right),
\label{eq1.3}
\end{equation}
and generalized $(PRS)_n$ to manifolds such that
\begin{equation}
\nabla_k P_{jl}  = 2A_k P_{jl}  + A_j P_{kl}  + A_l P_{kj}.
\label{eq1.4}
\end{equation}
The manifold is called {\em pseudo projective Ricci symmetric} 
and denoted by $(PWRS)_n$  \cite{[8]}.  
Recently another generalization of a $(PRS)_n$ was considered 
in \cite{[5]} and \cite{[10]}, whose Ricci tensor satisfies 
the condition 
\begin{equation}
\nabla_k R_{jl}  = (A_k  + B_k )R_{jl}  + A_j R_{kl}  + A_l R_{kj}, 
\label{eq1.5}
\end{equation}
The manifold is called 
{\em almost pseudo Ricci symmetric} and denoted by $A(PRS)_n$.  
In ref.\cite{[10]} the properties of conformally flat $A(PRS)_n$ were 
studied, pointing out their importance %of Pseudo Ricci symmetric manifolds 
in the theory of General Relativity. 

It seems worthwhile to introduce and study a new manifold structure 
that includes $(WRS)_n$, $(PRS)_n$ and $(PWRS)_n $ as special cases.\\

\begin{definition} 
A (0,2) symmetric tensor is a {\em generalized $Z$ tensor} if
\begin{equation}
Z_{kl}  = R_{kl}  + \phi\; g_{kl} \label{eq1.6},
\end{equation}
where $\phi $ is an arbitrary scalar function. The $Z$ scalar is
$Z=g^{kl}Z_{kl}= R+n\phi$.
\end{definition}
The classical $Z$ tensor is obtained with the choice $\phi = - \frac{1}{n} R$.
Hereafter we refer to the generalized $Z$ tensor simply as the $Z$ tensor.\\ 
The $Z$ tensor allows us to reinterpret several 
well known structures on Riemannian manifolds.\\
1) If $Z_{kl}=0$ the ($Z$-flat) manifold is an Einstein space, 
$R_{ij}=(R/n)g_{ij}$ \cite{[3]}.\\
2) If  $\nabla_i Z_{kl}= \lambda_i Z_{kl}$, the ($Z$-recurrent) manifold
is a generalized Ricci recurrent manifold \cite{[9a],[21]}: the
condition is equivalent to $\nabla_iR_{kl}=\lambda_iR_{kl}+(n-1)\,\mu_i\, 
g_{kl}$ where $(n - 1)\mu_i\equiv (\lambda_i -\nabla_i)\phi $. 
If moreover $0=(\lambda_i -\nabla_i)\phi$, a Ricci Recurrent manifold is 
recovered.\\
3) If $\nabla_k Z_{jl}  = \nabla_j Z_{kl}$ (i.e.  $Z$ is a Codazzi tensor, 
\cite{[13a]}) then
$\nabla_k R_{jl}-\nabla_jR_{kl}=(g_{kl}\nabla_j -g_{jl}\nabla_k)\phi$.
By transvecting with $g^{jl} $ we get 
$\nabla_k [R+2(n-1)\phi] = 0$ and, finally,
$$\nabla_k R_{jl}  - \nabla_j R_{kl}  = \frac{1}{2(n - 1)}
(g_{jl}\nabla_k - g_{kl}\nabla_j)R.$$
This condition defines a {\em nearly conformally symmetric} manifold, 
$(NCS)_n$. The condition was introduced and studied by Roter \cite{[23a]}. 
Conversely a $(NCS)_n$ has a Codazzi $Z$ tensor if 
$\nabla_k [R + 2(n - 1)\phi] = 0$.\\
4) Einstein's equations \cite{[14a]} with cosmological constant $\Lambda $ and 
energy-stress tensor $T_{kl}$ may be written as $Z_{kl}  = k T_{kl}$, 
where $\phi = - \frac{1}{2}R + \Lambda$, and $k$ is the 
gravitational constant.  
The $Z$ tensor may be thought of as a generalized Einstein 
gravitational tensor with arbitrary scalar function $\phi $.\\
Conditions on the energy-momentum tensor determine constraints on the
$Z$ tensor: the vacuum solution $Z = 0$ determines an Einstein space
with $\Lambda  = \frac{n - 2}{2n}\,R$;
conservation of total energy-momentum ($\nabla^l T_{kl}= 0$) gives
$\nabla^l Z_{kl}= 0$ and $\nabla_k (\frac{1}{2}R + \phi) = 0$;
the condition $\nabla_i Z_{kl}= 0$ describes a space-time with conserved 
energy-momentum density.\\

Several cases accomodate in a new kind of Riemannian manifold: 
\begin{definition} A manifold is {\em Weakly Z symmetric}, and 
denoted by $(WZS)_n$, if the generalized $Z$ tensor 
satisfies the condition:
\begin{equation}
\nabla_k Z_{jl}  = A_k Z_{jl}  + B_j Z_{kl}  + D_l Z_{kj} .\label{eq1.10}
\end{equation}
\end{definition}
If $\phi=0$ we recover a $(WRS)_n$ and its particular case 
$(PRS)_n$. If $\phi = -R/n$ (classical $Z$ tensor)
and if $A_k$ is replaced by $2A_k $, $B_k=D_k=A_k$, 
then $Z_{jl}=\frac{n - 1}{n}P_{jl}$ and the space reduces to a $(PWRS)_n$.\\

In sect.2 we obtain general properties of $(WZS)_n$
that descend directly from the definition and strongly depend on $Z_{ij}$ 
being singular or not. The two cases are examined in sections 3 and 4.
In sect.3 we study $(WZS)_n$ that are conformally or pseudo conformally
harmonic with $B-D\neq 0$; we show that $B-D$, after normalization,  
is a proper concircular vector. Sect.4 is devoted to $\mathrm{(WZS)}_n$ with 
non-singular Z tensor, and gives conditions for the closedness of the 1-form 
$A-B$ that involve various generalized curvature tensors.
In sect.5 we study conformally harmonic $\mathrm{(WZS)}_n$ and obtain the
explicit form of the Ricci tensor. In the conformally flat case we
also give the local form of the metric. 
 
\section{General properties} 
From the definition of a $(WZS)_n$ and its symmetries we obtain 
\begin{eqnarray}
0=\eta_j Z_{kl} - \eta_l Z_{kj}, \label{eq2.4}\\
\nabla_k Z_{jl}  - \nabla_j Z_{kl}  = \omega_k Z_{jl}  - \omega_j Z_{kl},
\label{eq2.2}
\end{eqnarray} 
with covectors $\omega_k=A_k-B_k$ and $\eta_k=B_k-D_k$ that will be 
used throughout.\\
Let's consider eq.(\ref{eq2.4}) first, it implies the following statements:
\begin{prop}\label{prop2.1}
In a $(WZS)_n$, if the $Z$ tensor is non-singular then $\eta_k=0$.
\begin{proof} If the $Z$ tensor is non singular, there exists a (2,0)
tensor $Z^{-1}$ such that $(Z^{-1})^{kh} Z_{kl}=\delta^h{}_l $. 
By transvecting eq.(\ref{eq2.4}) with $(Z^{-1})^{kh}$ we obtain
$\eta_j\delta_l{}^h  = \eta_l \delta_j{}^h$; put $h=l$ and sum to 
obtain $(n-1)\eta_j=0$.
\end{proof}
\end{prop}

\begin{prop}\label{prop2.2}
If $\eta_k\ne 0$ and the scalar $Z\neq 0$, then the $Z$ tensor 
has rank one:
\begin{equation}
Z_{ij} = Z\;\frac{\eta_i\eta_j}{\eta^k\eta_k} \label{rank1}
\end{equation}
\begin{proof}
Multiply eq.(\ref{eq2.4}) by $\eta^j$ and sum: $\eta^j\eta_j Z_{kl}
=\eta_l\eta^jZ_{kj}$. Multiply  eq.(\ref{eq2.4}) by $g^{jk}$ and sum:
$\eta^k Z_{kl}= Z \eta_l$. The two results imply the assertion.
\end{proof}
\end{prop}

The result translates to the Ricci tensor, whose expression is 
characteristic of {\em quasi Einstein} Riemannian manifolds \cite{[7]}, 
and generalizes the results of \cite{[11]}:
\begin{prop}\label{prop2.3}
A $(WZS)_n$ with $\eta_k\neq 0$, is a quasi Einstein manifold:
\begin{equation}
R_{ij}= -\phi\, g_{ij}+ Z\, T_i T_j, \qquad T_i=
\frac{\eta_i}{\|\eta \|},\label{Ricci}
\end{equation}
\end{prop} 

Next we consider eq.(\ref{eq2.2}). If $Z_{ij}$ is a Codazzi tensor, 
then the l.h.s. of the equation vanishes by definition, and the above 
discussion of 
eq.(\ref{eq2.4}) can be repeated. We merely state the result:
\begin{prop} \label{prop2.4}
In a $(WZS)_n$ with a Codazzi $Z$ tensor, if $Z$ is singular then 
$\omega_k\neq 0$. Conversely, if rank $[Z_{kl}]> 1$ then $\omega_k = 0$. 
\end{prop}

\section{Harmonic conformal or quasi conformal $\mathrm{(WZS)}_n$ with 
$\eta\ne 0$}
In this section we consider manifolds $(WZS)_n$ ($n > 3$) with $\eta_k\ne 0$,
and the property $\nabla_m C_{jkl}{}^m  = 0$ 
(i.e. harmonic conformal curvature tensor \cite{[3]}) or  
$\nabla_m W_{jkl}{}^m = 0$ (i.e. harmonic quasi conformal 
curvature tensor \cite{[28]}). We provide sufficient conditions for 
$\eta/\|\eta\|$ to be a proper concircular vector  \cite{[23],[26]}.\\

We begin with the case of harmonic conformal tensor. From the expression 
for the divergence of the conformal tensor,
\begin{equation}
\nabla_m C_{jkl}{}^m  = \frac{n - 3}{n - 2}
\left [\nabla_k R_{jl}  - \nabla_j R_{kl}
  + \frac{1}{2(n - 1)} (g_{kl}\nabla_j  - g_{jl}\nabla_k)R\right]\label{nablaC}
\end{equation}
we read the condition $\nabla_m C_{jkl}{}^m  = 0$:
\begin{equation}
 \nabla_k R_{jl}  - \nabla_j R_{kl}  = 
\frac{1}{2(n - 1)} (g_{jl}\nabla_k    - g_{kl}\nabla _j) R.
\label{eq3.2}
\end{equation}
We need the following theorem, whose proof given here 
is different from that in \cite{[12]} (see also \cite{[9]}):
\begin{thrm}\label{thrm3.1}
Let $M$ be a $n > 3$ dimensional manifold, with harmonic conformal
curvature tensor, and Ricci tensor $R_{kl}  = \alpha g_{kl}  + \beta T_k T_l $,
where $\alpha$, $\beta $ are scalars, and $T^kT_k=1$. If 
\begin{equation}
(T_j\nabla_k - T_k\nabla_j)\beta=0, \label{eq3.89}
\end{equation} 
then $T_k$ is a proper concircular vector.
\begin{proof}
Since $M$ is conformally harmonic, eq.(\ref{eq3.2}) gives:
\begin{equation}
 \beta[ \nabla_k (T_jT_l)  - \nabla_j (T_kT_l)]= 
\frac{1}{2(n - 1)} (g_{jl}\nabla_k    - g_{kl}\nabla_j )S,
\label{eq3.3}
\end{equation}
where $S = -(n-2)\alpha+\beta $, and condition (\ref{eq3.89})
was used. The proof is in four steps.\\
1) We show that $T^l\nabla_l T_k=0$: multiply eq.(\ref{eq3.3}) by 
$g^{jl}$ to obtain: a) $-\beta \nabla^l(T_kT_l)=\frac{1}{2}\nabla_k S $.
The result a) is multiplied by $T^k$ to give: b)
$-\beta \nabla_l T^l=\frac{1}{2}T^l\nabla_l S$. a) and b) combine to 
give: c) $-\beta T^l\nabla_l T_k = \frac{1}{2}[\nabla_k-T_kT^l\nabla_l]S$. 
Finally multiply eq.(\ref{eq3.3}) by $T^kT^l$ and use the property 
$T^l\nabla_kT_l=0$ to obtain:
$$\beta T^k\nabla_k T_j=\frac{1}{2(n - 1)}(T_jT^k\nabla_k-\nabla_j )S $$
which, compared to c) shows that d) $T^l\nabla_l T_k=0$ and 
$(T_jT^k\nabla_k-\nabla_j )S=0$.\\
2) We show that $T$ is a closed 1-form: multiply eq.(\ref{eq3.3}) by $T^l$
$$
\beta[ \nabla_k T_j  - \nabla_jT_k]= 
\frac{1}{2(n - 1)} (T_j\nabla_k - T_k\nabla_j)S. $$
$T$ is a closed form if the r.h.s. is null. This is proven by using identity a)
to write: $(T_j\nabla_k - T_k\nabla_j)S = 
-2\beta[ T_j\nabla^l(T_kT_l)- T_k\nabla^l(T_jT_l)]=0$ by property d).\\
3) With condition d) in mind, transvect eq.(\ref{eq3.3}) with $T^k$ and
obtain 
$$ -\beta\nabla_jT_l=\frac{1}{2(n-1)}(g_{jl}T^k\nabla_k - 
T_l\nabla_j)S $$
Use d) to replace $T_l\nabla_jS$ with 
$T_lT_jT^k\nabla_kS$. Then:
\begin{equation}
 \nabla_jT_l= f \,(T_jT_l-g_{jl}), \quad f\equiv
\frac{T^k\nabla_kS}{2\beta(n-1)} 
\end{equation}
which means that $T_k$ is a concircular vector. \\
4) We prove that $T_k$ is a proper concircular vector, i.e. $fT_k$  
is a closed 1-form: 
from d) by a covariant derivative we obtain
$\nabla_j\nabla_k S = (\nabla_j T_k )(T^l\nabla_l S)  
+ T_k \nabla_j (T^l\nabla_l S)$; subtract same equation with 
indices $k$ and $j$ exchanged. Since $T_k$ is a
closed 1-form we obtain: $T_k \nabla_j(T^l\nabla_l S) = 
T_j\nabla_k (T^l \nabla_l S)$. Multiply by $T^k$:
\begin{equation}
(T_j T^k\nabla_k -\nabla_j)(T^l \nabla_l S) = 0 \nonumber
\end{equation}
From the relation (\ref{eq3.89}), one obtains: 
$(T_kT^l\nabla_l-\nabla_k)\beta=0$. It follows that the scalar function $f$
has the property $\nabla_j f = \mu T_j$ where $\mu $ is a scalar 
function. Then the 1-form $f T_k$ is closed.
\end{proof}
\end{thrm}

With the identifications $\alpha =-\phi $ and $\beta =Z$, 
$T_i=\eta_i/\|\eta\|$ (see Prop. \ref{prop2.3}) 
the condition (\ref{eq3.89}) is $(\eta_j \nabla_k - \eta_k\nabla_j) Z=0$. 
Since $Z =S -(n-2)\phi $ 
and $(\eta_j\nabla _k -\eta_k\nabla_j) S=0$, the condition can be 
rewritten as $(\eta_j\nabla _k -\eta_k \nabla _j) \phi =0$.
Thus we can state the following:

\begin{thrm}\label{thrm3.2}
In a $(WZS)_n$ manifold with $\eta_k\neq 0$ and harmonic conformal 
curvature tensor, if 
\begin{equation}
(\eta_j \nabla_k - \eta_k \nabla_j )\phi=0 \label{cond}
\end{equation}
then $\eta_i/\| \eta \|$ is a proper concircular vector.
\end{thrm}

\begin{rem}\label{rem3.3} 
If $\phi=0$ or $\nabla_k \phi=0$, 
the condition (\ref{cond}) is fulfilled automatically. 
In the case $\phi=0$ we recover a $(WRS)_n$
manifold (and the results of refs \cite{[9],[11]}).
\end{rem}

Now we consider the case of a $(WZS)_n$ manifold with harmonic quasi 
conformal curvature tensor.
In 1968 Yano and Sawaki \cite{[28]} defined and studied a tensor $W_{jkl}{}^m$
on a Riemannian manifold of dimension $n>3$, 
which includes as particular cases the conformal 
curvature tensor $C_{jkl}{}^m$, eq.(\ref{eq1.2}), and the concircular 
curvature tensor
\begin{equation}
\tilde C_{jkl}{}^m  = R_{jkl}{}^m  + \frac{R}{n(n-1)}(\delta_j{}^m g_{kl}  
- \delta_k{}^m g_{jl} ).\label{eq3.17}
\end{equation}
The tensor is known as the {\em quasi conformal} curvature tensor:
\begin{equation}
W_{jkl}{}^m  = -(n - 2)\,b\,C_{jkl}{}^m  + [a + (n - 2)b]\tilde C_{jkl}{}^m;
\label{eq3.16}
\end{equation}
$a$ and $b$ are nonzero constants. From the expressions (\ref{nablaC})
and (\ref{nablaCtilde}) we evaluate
\begin{equation}
\nabla_m W_{jkl}{}^m  = (a + b)\nabla_m R_{jkl}{}^m  + 
\frac{2a-b(n - 1)(n - 4)}{2n(n - 1)}
(g_{kl}\nabla_j  - g_{jl}\nabla _k) R.\label{eq3.18}
\end{equation}
A manifold is {\em quasi conformally harmonic} 
if $\nabla_m W_{jkl}{}^m =0$.
By transvecting the condition with $g^{jk}$ we get:
\begin{equation}
(1 - 2/n)[a + b(n - 2)]\;\nabla _j R \,= \,0,\label{condQCH} 
\end{equation}
which means that either $a + b(n - 2) = 0$ or $\nabla_j R = 0$.
The first condition implies $W=C$, and gives back the harmonic 
conformal case. If $\nabla_j R = 0$ it is $\nabla_m R_{jkl}{}^m =0$
by (\ref{eq3.18}), and the equations in the proof of theorem \ref{thrm3.1} 
simplify and we can state the following 
(analogous to theorem \ref{thrm3.2}):

\begin{thrm}\label{thrm3.3}
Let $(WZS)_n$ be a quasi conformally harmonic manifold, with 
$\eta_k\neq 0$. If $(\eta_j\nabla_k -\eta_k\nabla_j)\phi =0 $, 
then $\eta/\|\eta\|$ is a proper concircular vector.
\end{thrm}

\section{$\mathrm{(WZS)}_n$ with non-singular Z tensor: 
conditions for closed $\omega $} 
In this section we investigate in a $(WZS)_n$ $(n > 3)$ the conditions 
the 1-form $\omega_k$ to be closed: $\nabla_i\omega_j-\nabla_j\omega_i=0$.
We need: 
\begin{lem}[Lovelock's differential identity, \cite{[19],[20]}]\label{lem4.1} 
In a Riemannian manifold the following identity is true:
\begin{eqnarray}
\nabla_i\nabla_m R_{jkl}{}^m  + \nabla_j\nabla_m R_{kil}{}^m  + 
\nabla_k\nabla_m R_{ijl}{}^m \nonumber\\ 
= -R_{im} R_{jkl}{}^m  - R_{jm} R_{kil}{}^m  
- R_{km} R_{ijl}{}^m \label{eq4.1}
\end{eqnarray}
\end{lem}
\noindent
and also the contracted second Bianchi identity in the form
\begin{equation}
\nabla_m R_{jkl}{}^m  = \nabla_k Z_{jl}  - \nabla_j Z_{kl}  
+ (g_{kl}\,\nabla_j  -  g_{jl}\,\nabla _k) \phi.\label{eq4.2}
\end{equation}
Now we prove the relevant theorem (see also \cite{[20]}):
\begin{thrm}\label{thrm4.1}
In a $(WZS)_n$ ($n>3$) with non singular $Z$ tensor, $\omega_k$ is a 
closed 1-form if and only if:
\begin{equation}
R_{im} R_{jkl}{}^m  + R_{jm} R_{kil}{}^m  + R_{km} R_{ijl}{}^m  = 0.
\label{eq4.4}
\end{equation}
\begin{proof} 
The covariant derivative of eq.(\ref{eq4.2}) and eq.(\ref{eq2.2}) give:
$\nabla_i \nabla_m R_{jkl}{}^m  = (\nabla_i\omega_k)Z_{jl}  
+ \omega_k (\nabla_i Z_{jl}) - (\nabla_i\omega_j )Z_{kl}  
- \omega_j (\nabla_i Z_{kl} ) + ( g_{kl}\nabla_i \nabla_j \phi  
-  g_{jl}\nabla _i \nabla _k \phi ). $
Cyclic permutations of the indices $i,j,k$ are made, and the 
resulting three equations are added:
\begin{eqnarray}
&&\nabla_i\nabla_m R_{jkl}{}^m  + \nabla_j\nabla_m R_{kil}{}^m  + 
\nabla_k\nabla_m R_{ijl}{}^m   \nonumber\\
&&=(\nabla_i\omega_k  - \nabla_k\omega_i) Z_{jl}  
+ (\nabla_j\omega_i - \nabla_i\omega_j)Z_{kl}  
+ (\nabla_k\omega_j - \nabla_j\omega_k) Z_{il}\nonumber \\
&&\quad + \omega_j (\nabla_k Z_{il}  - \nabla_i Z_{kl}) + 
\omega_k (\nabla_i Z_{jl}  - \nabla_j Z_{il} ) + 
\omega_i (\nabla_j Z_{kl}  - \nabla_k Z_{jl} ).\nonumber
\end{eqnarray}
Cancellations occur by eq.(\ref{eq2.2}). By lemma \ref{lem4.1}, 
one obtains:
\begin{eqnarray}
&&-R_{im} R_{jkl}{}^m  - R_{jm} R_{kil}{}^m  - R_{km} R_{ijl}{}^m \nonumber \\
&&=(\nabla_i\omega_k - \nabla_k\omega_i)Z_{jl} + 
(\nabla_j\omega_i - \nabla_i\omega_j)Z_{kl} + 
(\nabla_k\omega_j - \nabla_j\omega_k)Z_{il}.\nonumber
\end{eqnarray}
If $\omega_k $ is a closed 1-form then eq.(\ref{eq4.4}) is fulfilled.
Conversely, suppose that eq.(\ref{eq4.4}) holds: if the $Z$ tensor is 
non singular, there is a $(2,0)$ tensor such that 
$Z_{kl}(Z^{-1})^{km}= \delta_l{}^m $. Multiply the last equation by 
$(Z^{ - 1} )^{hl}$:
$(\nabla_i\omega_k - \nabla_k\omega_i)\delta_j{}^h + 
(\nabla_j\omega_i - \nabla_i\omega_j)\delta_k{}^h + 
(\nabla_k\omega_j - \nabla_j\omega_k)\delta_i{}^h  = 0$.
Set $h=j$ and sum: 
$(n - 2)(\nabla_i\omega_k - \nabla_k\omega_i) = 0$. Since $n > 2$, 
$\omega_k$ is a closed 1-form.
\end{proof}
\end{thrm}

\begin{rem}
By Lovelock's identity, the condition (\ref{eq4.4}) is obviously true
if $\nabla_m R_{ijk}{}^m=0$, i.e. the $(WZS)_n$ is a harmonic manifold.
However, we have shown in ref.\cite{[20]} that there is a broad class of 
generalized curvature tensors for which the case $\nabla_m K_{ijk}{}^m=0$ 
implies the same condition. This class includes several well known
curvature tensors, and is the main subject of this section.
\end{rem}

\begin{definition}
A tensor $K_{jkl}{}^m$ is a {\em generalized curvature tensor}\footnote{The
notion was introduced by Kobayashi and Nomizu \cite{[17a]}, but with the 
further antisymmetry in the last pair of indices.} if:\\
1) $K_{jkl}{}^m  =  - K_{kjl}{}^m$,\\ 
2) $K_{jkl}{}^m  + K_{klj}{}^m  + K_{ljk}{}^m  = 0$.
\end{definition}
The second Bianchi identity does not hold in general, and 
is modified by a tensor source $B_{ijkl}{}^m$ that depends on the 
specific form of the curvature tensor:  
\begin{equation}
\nabla_i K_{jkl}{}^m  + \nabla_j K_{kil}{}^m  
+ \nabla_k K_{ijl}{}^m  = B_{ijkl}{}^m  \label{Btensor}
\end{equation}
\begin{prop}[\cite{[20]}]\label{thrm4.2}
If $K_{jkl}{}^m$ is a generalized curvature tensor such that
\begin{eqnarray}
\nabla _m K_{jkl}{}^m  = A\nabla_m R_{jkl}{}^m  + B(a_{lk}\nabla_j - a_{lj}
\nabla_k) \psi ,\label{eq4.11}
\end{eqnarray}
where $A\neq 0,\,B$ are constants, $\psi$ is a scalar field,
and $a_{ij}$ is a symmetric $(0,2)$ Codazzi tensor (i.e. 
$\nabla_ia_{kl}= \nabla_ka_{il}$), then the following relation holds:
\begin{eqnarray}
&&\nabla_i\nabla_m K_{jkl}{}^m  + \nabla_j\nabla_m K_{kil}{}^m  + 
\nabla_k\nabla_m K_{ijl}{}^m\label{eq4.12} \\  
&& \qquad= -A(R_{im} R_{jkl}{}^m  
+ R_{jm} R_{kil}{}^m  + R_{km} R_{ijl}{}^m ).\nonumber
\end{eqnarray}
\end{prop}
\begin{rem}
%The existence of non-trivial Codazzi tensors has important geometrical 
%and topological consequences \cite{[3a]}. 
In \cite{[13a]} it is proven that any smooth manifold carries a metric  
such that $(M,g)$ admits a non trivial Codazzi tensor (i.e.
proportional to the metric tensor) and 
the deep consequences on the structure of the curvature operator are 
presented (see also \cite{[20a]}).\\
Given a Codazzi tensor it is possible to exhibit a $K$ tensor that 
satisfies the condition (\ref{eq4.11}):
\begin{equation}
     K_{jkl}{}^m = A\;R_{jkl}{}^m + B\,\psi\, 
(\delta_j{}^m a_{kl}-\delta_k{}^m a_{jl} ).\label{eq4.13}
\end{equation} 
Its trace is: $K_{kl}= - K_{mkl}{}^m = A\,R_{kl} - B (n-1)\psi\,a_{kl}$. 
Note that for $a_{kl}=g_{kl}$ the tensor $K_{kl}$ is up to a factor 
a $Z$ tensor. Thus $Z$ tensors arise naturally from the invariance 
of Lovelock's identity.
\end{rem}

\begin{rem} In the literature one meets generalized curvature tensors whose 
divergence has the form (\ref{eq4.11}), with trivial Codazzi tensor:
\begin{equation}
\nabla_m K_{jkl}{}^m = A\,\nabla_m R_{jkl}{}^m + B (g_{kl}\nabla_j  
-g_{jl}\nabla_k) R. \label{eq4.15}
\end{equation}
They are the projective curvature tensor $P_{jkl}{}^m$ \cite{[15]}, 
the conformal curvature tensor $C_{jkl}{}^m$ \cite{[22]}, the concircular 
tensor $\tilde C_{jkl}{}^m$ \cite{[23],[26]}, the conharmonic tensor 
$N_{jkl}{}^m$ \cite{[21],[24]} and the quasi conformal tensor 
$W_{jkl}{}^m$ \cite{[28]}.
\end{rem}

\begin{definition}
A manifold is $K$-harmonic if $\nabla_m K_{jkl}{}^m =0$.
\end{definition}

\begin{prop}
In a $K$-harmonic manifold, if $K$ is of type (\ref{eq4.15}) and
$A\neq 2(n-1)B$, then $\nabla_j R=0$.
\begin{proof}
By transvecting eq.(\ref{eq4.15}) with $g^{kl}$ and by the second 
contracted Bianchi identity, we obtain $\frac{1}{2}[A -2(n-1)B]\nabla_j R=0$.
\end{proof}
\end{prop}

Hereafter, we specialize to $(WZS)_n$ manifolds with non singular $Z$ tensor, 
and with a generalized curvature tensor of the type (\ref{eq4.15}). From 
eqs. (\ref{eq4.2}) and (\ref{eq2.2}) we obtain:
\begin{equation}
\nabla_m K_{jkl}{}^m = A(\omega_k Z_{jl} - \omega_j Z_{kl}) + (g_{kl}\nabla_j
- g_{jl}\nabla_k)(A\phi + B\,R). \label{eq4.17}
\end{equation}
Then, the manifold is $K$-harmonic if:
\begin{equation}
A(\omega_k Z_{jl} - \omega_j Z_{kl}) = ( g_{jl}\nabla_k 
 - g_{kl}\nabla_j)\, (A\phi + B\,R ). \label{eq4.18}
\end{equation}

\begin{lem}\label{prop4.2}
In a $K$-harmonic $(WZS)_n$ with non singular $Z$ tensor:\\ 
1) $\omega_k=0$ if and only if $\nabla_k(A\phi +B R ) = 0$;\\
2) If $A\neq 2(n-1)B$, then $\omega_k=0$ if and only if $\nabla_k\phi =0$.
\begin{proof}
If $\nabla_k(A\phi + BR) = 0$ then $\omega_k Z_{jl} = \omega_j Z_{kl}$:
if the $Z$ tensor is non singular,
by transvecting with $(Z^{-1})^{lh}$ we obtain
$\omega_j\delta^h{}_k  = \omega_k \delta^h{}_j$. Now put $h=j$ and sum 
to obtain $(n-1)\omega_k=0$. On the other hand if $\omega_k=0$ 
eq.(\ref{eq4.18}) gives
$[ g_{jl}\nabla_k - g_{kl}\nabla_j] (A\phi + B R )=0$ 
and transvecting with $g^{kl}$ we get the result.\\
If $A\neq 2B(n-1)$ then $\nabla_k R = 0$ and part 1) applies.   
\end{proof}
\end{lem}

\begin{thrm}\label{thrm4.3b}
In a $K$-harmonic $(WZS)_n$ with non-singular $Z$ tensor  and $K$ of type
(\ref{eq4.15}), if $\omega\neq 0$ then $\omega $ is a closed 1-form.
\end{thrm}

This theorem extends theorem \ref{thrm4.1} (where $K=R$), 
and has interesting corollaries according to the various choices 
$K=C,\,W,\,P,\,\tilde C,\,N$.

\begin{cor}\label{cor99}
Let $(WZS)_n$ have non singular $Z$ tensor and $\omega\neq 0$. 
If $\nabla_m K_{jkl}{}^m = 0$, and $K= P,\,\tilde C,\,N$, 
then $\omega $ is a closed 1-form.
\begin{proof}
1) Harmonic conformal curvature: $\nabla_m C_{jkl}{}^ m = 0$. Note that
in this case $A=2B(n-1)$; theorem \ref{thrm4.3b} applies.\\
2) Harmonic quasi conformal curvature: $\nabla_m W_{jkl}{}^m=0$: 
Eq.(\ref{condQCH}) gives either $\nabla_j R=0$ or $a + b(n-2) = 0$.
If $\nabla_j R = 0$ then $\nabla_m R_{jkl}{}^m = 0$ and 
theorem \ref{thrm4.1}. If $a + b(n-2) = 0$ it is 
$\nabla_m C_{jkl}{}^m = 0$ and case 1) applies.\\
3) Harmonic projective curvature: $\nabla_m P_{jkl}{}^m = 0$.
The components of the projective curvature tensor are \cite{[15],[24]}:
\begin{equation*}
P_{jkl}{}^m  = R_{jkl}{}^m  + \frac{1}{n - 1}(\delta_j{}^m R_{kl}  
- \delta_k{}^m R_{jl} ).
\end{equation*}
One evaluates $\nabla_m P_{jkl}{}^m  = \frac{n - 2}{n - 1}\nabla_m R_{jkl}{}^m$,
and theorem \ref{thrm4.1} applies.\\
4) Harmonic concircular curvature: $\nabla_m\tilde C_{jkl}{}^m  = 0$.
The concircular curvature tensor is given in eq.(\ref{eq3.17}), 
\cite{[23],[26]}. Its divergence is
\begin{equation}
 \nabla _m \tilde C_{jkl}{}^m  = \nabla _m R_{jkl} ^m  + \frac{1}
{n(n - 1)}(g_{kl}\nabla_j  - g_{jl}\nabla_k)R\label{nablaCtilde}
\end{equation}
Theorem \ref{thrm4.3b} applies.\\
5) Harmonic conharmonic curvature: $\nabla_m N_{jkl}{}^m  = 0$. 
The conharmonic curvature tensor \cite{[21],[24]} is:
\begin{equation*}
N_{jkl}{}^m  = R_{jkl}{}^m  + \frac{1}{n - 2}(\delta_j{}^m R_{kl}  
- \delta_k{}^m R_{jl}  + R_j^m g_{kl}  - R_k^m g_{jl} ).
\end{equation*}
A covariant derivative and the second contracted Bianchi identity give:
$$\nabla_m N_{jkl}{}^m  = \frac{n - 3}{n - 2}\nabla_m R_{jkl}{}^m  + 
\frac{1}{2(n - 2)}(g_{kl}\nabla_j - g_{jl}\nabla _k)R. $$
Theorems \ref{thrm4.3b} applies.
\end{proof}
\end{cor}

There are other cases where the 1-form $\omega_k$ is closed for a 
$(WZS)_n$ manifold.

\begin{definition}[\cite{[20],[17]}]
A $n$-dimensional Riemannian manifold is {\em $K$-recurrent}, $(KR)_n$,
if the generalized curvature tensor is recurrent, 
$\nabla_i K_{jkl}{}^m  = \lambda_i K_{jkl}{}^m $,
for some non zero covector $\lambda_i$. 
\end{definition}
\begin{thrm}[\cite{[20]}]\label{thrm4.4}
In a $(KR)_n$, if $\lambda_i$ is closed then:
\begin{equation}
R_{im} R_{jkl}{}^m  + R_{jm} R_{kil} {}^m  + R_{km} R_{ijl}{}^m  = 
-\frac{1}{A}\nabla_m B_{ijkl}{}^m .\label{eq4.25}
\end{equation}
where $B$ is the source tensor in eq.(\ref{Btensor}). 
In particular, for $K = C,\,P,\,\tilde C,\,N,\,W$ the tensor 
$\nabla_m B_{ijkl}{}^m$ either vanishes or is proportional to the l.h.s.
of eq.(\ref{eq4.25}).
\end{thrm}

\begin{cor}\label{cor4.6}
Let $(WZS)_n$ have non singular $Z$ tensor, and be $K$ recurrent with 
closed $\lambda_i$. If $K=C,\,P,\,\tilde C,\,N,\,W$, then $\omega $ is 
a closed 1-form.
\end{cor}

\begin{definition}
A Riemannian manifold is {\em pseudosymmetric in the sense of 
R.~Deszcz} \cite{[14]} if the following condition holds:
\begin{eqnarray}
(\nabla_s\nabla_i  - \nabla_i\nabla_s)R_{jklm} = L_R\, (g_{js} R_{iklm}
- g_{ji} R_{sklm}  + g_{ks} R_{jilm}  - g_{ki} R_{jslm}   \nonumber\\
+ g_{ls} R_{jkim}  - g_{li} R_{jksm}  + g_{ms} R_{jkli} 
- g_{mi} R_{jkls} ), \label{eq4.26}
\end{eqnarray} 
where $L_R $ is a non null scalar function.
\end{definition}
In ref.\cite{[20]} the following theorem is proven:
\begin{thrm}\label{thrm4.7}
In a Riemannian manifold which is pseudosymmetric in 
the sense of R.~Deszcz, it is
$ R_{im} R_{jkl}{}^m  + R_{jm} R_{kil}{}^m  + R_{km} R_{ijl}{}^m  = 0.$
\end{thrm}
Then we can state the following:
\begin{prop}
In a $(WZS)_n$ which is pseudosymmetric in the sense of R.~Deszcz, 
if the $Z$ tensor is non-singular then $\omega_k$ is a closed 1-form.
\end{prop}

\begin{definition} A Riemannian manifold is 
{\em generalized Ricci pseudosymmetric in the sense of R.~Deszcz}, \cite{[13]},
if the following condition holds:
\begin{eqnarray}
(\nabla_s\nabla_i - \nabla_i\nabla_s) R_{jklm}  = L_S (R_{js} R_{iklm}  - 
R_{ji} R_{sklm}  + R_{ks} R_{jilm}  - R_{ki} R_{jslm}  +  \nonumber\\
+ R_{ls}R_{jkim}-R_{li}R_{jksm}+R_{ms}R_{jkli}-R_{mi}R_{jkls}), \label{eq4.27}
\end{eqnarray}
where $L_S$ is a non null scalar function. 
\end{definition}

\begin{thrm}\label{thrm4.9}
In a generalized Ricci pseudosymmetric manifold in the sense of 
R.~Deszcz, it is either 
$L_S=-\frac{1}{3}$, or $R_{im} R_{jkl}{}^m  + R_{jm} R_{kil}{}^m  + 
R_{km} R_{ijl}{}^m  = 0.$
\begin{proof}
Equation (\ref{eq4.27}) is transvected with $g^{mj} $ to obtain
$$ (\nabla_s\nabla_i - \nabla_i\nabla_s)R_{kl} = 
L_S [R_{im}(R_{skl}{}^m+R_{slk}{}^m ) - R_{sm} (R_{ikl}{}^m  + R_{ilk}{}^m )].$$
Then:
\begin{eqnarray}
&&(\nabla_i\nabla_k -\nabla_k\nabla_i)R_{jl}  
+ (\nabla_j\nabla_i-\nabla_i\nabla_j )R_{kl}  
+ (\nabla_k\nabla_j-\nabla_j\nabla_k )R_{il}\nonumber \\
&& = 3L_S (R_{im} R_{jkl}{}^m  + R_{jm} R_{kil}{}^m  + R_{km} R_{ijl}{}^m )
\nonumber
\end{eqnarray} 
By Lovelock's identity (\ref{lem4.1}), 
the l.h.s. of the previous equation is:
\begin{eqnarray}
&&\nabla_i\nabla_m R_{jkl}{}^m  + \nabla_j\nabla_m R_{kil}{}^m  
+ \nabla_k\nabla_m R_{ijl}{}^m \nonumber\\
&&= (\nabla_i\nabla_k -\nabla_k\nabla_i)R_{jl}  
+ (\nabla_j\nabla_i-\nabla_i\nabla_j )R_{kl}  
+ (\nabla_k\nabla_j-\nabla_j\nabla_k )R_{il}\nonumber \\
&&= -R_{im} R_{jkl}{}^m  - R_{jm} R_{kil}{}^m  - R_{km} R_{ijl}{}^m . \nonumber 
\end{eqnarray} 
Compare the two results and conclude that either
$L_S =-\frac{1}{3}$, or $R_{im} R_{jkl}{}^m  + R_{jm}R_{kil}{}^m +
R_{km}R_{ijl}{}^m  = 0.$ 
\end{proof}
\end{thrm}

Finally we state:
\begin{prop}
In a $(WZS)_n$ which is also a generalized Ricci pseudosymmetric manifold
in the sense of R.Deszcz, if the $Z$ 
tensor is non-singular and $L_S \ne -\frac{1}{3}$, then $\omega _k $
is a closed 1-form.
\end{prop}

\section{Conformally harmonic $\mathrm{(WZS)}_n$: form of the
Ricci tensor}

In this section we study conformally harmonic $(WZS)_n$ in depth. 
We show the existence of a proper concircular vector 
in such manifolds, and obtain the form of the Ricci tensor. The proof 
only requires the $Z$ tensor to be non singular. For the conformally flat 
case, in particular, we give the explicit local form of the metric tensor.

The condition $\nabla_m C_{jkl}{}^m = 0$ is eq.(\ref{eq3.2}) which, by using
$R_{ij}=Z_{ij}-g_{ij}\;\phi$ and the property eq.(\ref{eq2.2}), becomes:
\begin{equation}
\omega_k Z_{jl}  - \omega_j Z_{kl}  = \frac{1}{2(n - 1)}
(g_{jl}\nabla_k  - g_{kl}\nabla_j )[R + 2(n - 1)\phi ].
\label{eq5.1}
\end{equation}
This is the starting point for the proofs. By prop \ref{prop4.2}, since
$Z$ is non singular,  
$\omega_k\neq 0$ if and only if $\nabla_k [R + 2(n - 1)\phi ]\neq 0$.

\begin{rem}\label{rem5.1}
1) The condition $\nabla_m C_{jkl}{}^m = 0$ implies that the 
manifold is a $(NCS)_n$.\\ 
2) If $\nabla_k [R+2(n-1)\phi]=0$ the $Z$ tensor is a Codazzi tensor. 
\end{rem}

The following theorem generalizes a result in \cite{[10]} for $A(PRS)_n$:
\begin{thrm}\label{thrm5.1}
In a conformally harmonic $(WZS)_n$ the 1-form $\omega $ is an eigenvector 
of the $Z$ tensor. 
\begin{proof}
By transvecting eq.(\ref{eq5.1}) with $g^{kl} $ we obtain
\begin{equation}
 \omega_j Z - \omega^m Z_{jm}  = \frac{1}{2}\nabla _j [R + 2(n - 1)\phi ];
\label{eq5.2}
\end{equation}
the result is inserted back in eq.(\ref{eq5.1}),
\begin{equation}
\omega_k Z_{jl}  - \omega_j Z_{kl}  = \frac{1}{(n - 1)}
[(\omega_k Z - \omega^m Z_{km})g_{jl}- (\omega_j Z - \omega^m Z_{jm})g_{kl}],
\nonumber
\end{equation}
and transvected with $\omega^j\omega^l$ to obtain 
$\omega_k (\omega^j\omega^l Z_{jl})  = (\omega_j \omega^j)\omega ^l Z_{kl}$.
The last equation can be rewritten as:
$ Z_{kl}\omega^l  = \zeta\omega_k$
\end{proof}
\end{thrm}

Now eq.(\ref{eq5.2}) simplifies:
$ \omega_j (\zeta - Z) =  - \frac{1}{2}\nabla_j [R + 2(n - 1)\phi ]$.
The result is a natural generalization of a similar one 
given in ref.\cite{[10]} for $A(PRS)_n$.

\begin{thrm}\label{thrm5.2}
Let $M$ be a conformally harmonic $(WZS)_n$. Then:\\
1) $M$ is a quasi Einstein manifold;\\
2) if the $Z$ tensor is non singular and if 
$(\omega_j\nabla_k -\omega_k\nabla_j)\phi= 0$, then:
\begin{equation}
(\omega_j\nabla_k-\omega_k\nabla_j)\left[\frac{n\zeta - Z}{n - 1}\right] = 0,
\label{eq5.55}
\end{equation}
and $M$ admits a proper concircular vector.
\begin{proof}
Eq.(\ref{eq5.1}) is transvected with $\omega^j$ and theorem \ref{thrm5.1}
is used to show that
$$ R_{kl}  = \left[\frac{Z - \zeta}{n - 1} - \phi\right] g_{kl} + 
\left[\frac{n\zeta - Z}{n - 1}\right]
\frac{\omega_k \omega_l}{\omega_j\omega^j},$$
i.e. $R_{kl}$ has the structure $\alpha g_{kl}  + \beta T_k T_l$ and the 
manifold is quasi Einstein \cite{[7]}. By transvecting 
eq.(\ref{eq4.2}) with $g^{jl}$ we obtain 
\begin{equation*}
\frac{1}{2}\nabla_k Z + \frac{n - 2}{2}\nabla_k \phi  = \omega_k Z - 
\omega^l Z_{kl}.
\end{equation*}
This and theorem (\ref{thrm5.1}) imply:
\begin{equation}
\frac{1}{2}\nabla_k Z + \frac{n - 2}{2}\nabla_k \phi  = \omega_k (Z - \zeta).
\label{eq5.10}
\end{equation}
A covariant derivative gives $
\frac{1}{2}\nabla_j \nabla_k Z + \frac{n - 2}{2}\nabla_j \nabla_k \phi  
= \nabla_j[ \omega_k (Z - \zeta )]$. 
Subtract the equation with indices $k$ and $j$ exchanged:
\begin{equation*}
(Z-\zeta )(\nabla_j\omega_k  - \nabla_k\omega_j )+ (\omega_k\nabla_j 
- \omega_j\nabla_k) (Z - \zeta) = 0.
\end{equation*}
According to corollary \ref{cor99}, in a conformally harmonic $(WZS)_n$ 
with non singular $Z$ the 1-form $\omega_k $ is closed. Then
\begin{equation}
  (\omega_k \nabla_j - \omega_j\nabla_k) (Z - \zeta ) = 0 \label{eq5.14}
\end{equation} 
Multiply eq.(\ref{eq5.10}) by $\omega_j$ and subtract from it the 
equation with indices $k$ and $j$ exchanged: 
$(\omega_j\nabla_k  - \omega_k\nabla_j )Z + (n - 2)(\omega_j\nabla_k  
- \omega_k\nabla_j)\phi = 0.$ Suppose that $\omega_k$, besides being a 
closed 1-form, has the property
$(\omega_j\nabla_k -\omega_k\nabla_j)\phi = 0$, then one obtains 
the further equation:
\begin{equation}
 (\omega_k \nabla_j - \omega_j\nabla_k) Z = 0. \label{eq5.15}
\end{equation} 
Eqs.(\ref{eq5.14},\ref{eq5.15}) imply the assertion eq.(\ref{eq5.55}).
The existence of a proper concircular vector follows from
Theorem \ref{thrm3.1}.
\end{proof}
\end{thrm}

Let us specialize to the case $C_{ijk}{}^m=0$ (conformally flat $(WZS)_n$).\\ 
It is well known \cite{[1]} that if a conformally flat space 
admits a proper concircular vector, then the space is subprojective 
in the sense of Kagan.\\
From theorem \ref{thrm5.2} we state the following:

\begin{thrm}\label{thrm5.4}
Let $(WZS)_n$ $(n > 3)$ be conformally flat with nonsingular $Z$ tensor 
and $(\omega_j\nabla_k -\omega_k\nabla_j)\phi  = 0$, then the manifold 
is a subprojective space.
\end{thrm}

In \cite{[27]} K. Yano proved that a necessary and sufficient condition 
for a Riemannian manifold to admit a concircular vector, is that there is 
a coordinate system in which the first fundamental form may be written as:
\begin{equation}
ds^2  = (dx^1 )^2  + e^{q(x^1 )} g^*_{\alpha \beta }(x^2,\ldots ,x^n) 
dx^\alpha  dx^\beta, \label{metric}
\end{equation}
where $\alpha ,\beta = 2,\ldots,n$. Since a conformally flat $(WZS)_n$ 
with non singular $Z$ tensor admits a proper 
concircular vector field, this space is the warped product $1 \times e^q M^*$,  
where $(M^*,g^*)$ is a $(n-1)-$dimensional Riemannian manifold. 
Gebarosky \cite{[16]} 
proved that the warped product $1 \times e^q M^*$ has the metric structure
(\ref{metric}) if and only if $M^*$ is Einstein. 
Thus the following theorem holds:

\begin{thrm}\label{5.5}
Let $M$ be a $n$ dimensional conformally flat $(WZS)_n$ $(n > 3)$.  
If $Z_{kl}$ is non singular and $(\omega_j \nabla_k -\omega_k \nabla_j)\phi=0$,
then $M$ is the warped product $1\times e^q M^*$, where  $M^*$ is Einstein.
\end{thrm}

\vfill
\end{document}